\documentclass[11pts,a4]{article}
\usepackage{amsmath,amssymb,amsfonts,amsthm,euscript,amsbsy,epic}
\usepackage[dvips]{graphicx}

\usepackage[usenames]{color}

\newcommand{\ch}[1]{{#1}}

\newcommand{\wh}[1]{\widehat{#1}}
\newcommand{\wt}[1]{\widetilde{#1}}
\newtheorem{THM}{Theorem}[section]

\newtheorem{REM}[THM]{Remark}
\newtheorem{LEM}[THM]{Lemma}

\newcommand{\R}[1]{\frac{1}{#1}}

\def\Rd{\mathbb{R}^d}

\def\ra{\rightarrow}
\def\t0{\wt{t}_0(0)}

\def \mP{\mathbb{P}}

\def\hlf{\frac{1}{2}}

\def\dtwo{\frac{d-2}{2}}
\def\dfour{\frac{d-4}{2}}

\def\deight{\frac{d-8}{2}}

\def\N{\mathbb{N}}
\def\Z{\mathbb{Z}}

\newcommand{\nn}{\nonumber}

\newcommand{\eqn}[1]{\begin{equation*}#1\end{equation*}}
\newcommand{\seqn}[1]{\begin{equation*}\begin{split}#1
  \end{split}\end{equation*}}
\newcommand{\eqnlab}[2]{\begin{equation}\label{#1}#2\end{equation}}
\newcommand{\seqnlab}[2]{\begin{equation}\label{#1}\begin{split}#2
  \end{split}\end{equation}}

\newcommand{\sm}[3]{\frac{#1}{#2^{#3}}}
\newcommand{\mc}[1]{\mathcal{#1}}

\title{An extension of the inductive approach to the lace expansion}
\author{Remco van der Hofstad
\thanks{Department of Mathematics and
Computer Science, Eindhoven University of Technology, P.O.\ Box
513, 5600 MB Eindhoven, The Netherlands.
E-mail: {rhofstad@win.tue.nl, holmes@eurandom.tue.nl}
}
\and Mark Holmes$\;^*$ \and Gordon Slade
\thanks{Department of Mathematics, University of
British Columbia, Vancouver, BC V6T 1Z2, Canada.
E-mail: {slade@math.ubc.ca}
}
}

\date{June 3, 2007}

\begin{document}
\maketitle
\begin{abstract}
We extend the inductive approach to the lace expansion, previously
developed to study models with critical dimension $4$, to be
applicable more generally.  In particular, the result of this note
has recently been used to prove Gaussian asymptotic behaviour for
the Fourier transform of the two-point function for sufficiently
spread-out lattice trees in dimensions $d>8$, and it is potentially also
applicable to percolation in dimensions $d>6$.
\end{abstract}

%%%%%%%%%%%%%%%%%%%%%%%%%%

\section{Motivation}
The lace expansion has been used since the mid-1980s to study a wide variety
of problems in high-dimensional probability, statistical mechanics, and
combinatorics \cite{Slad06}.  One of the most flexible approaches to the lace
expansion is the inductive method, first developed in \cite{HHS98} in the
context of weakly self-avoiding walks in dimensions $d>4$, and subsequently
extended to a much more general setting in \cite{HS02}.
The inductive approach of \cite{HS02} was successfully used
to prove Gaussian asymptotic behavior for the Fourier transform of the critical
two-point function $c_n(x;z_c)$ for a sufficiently spread-out model of self-avoiding
walk in dimensions $d>4$ \cite{HS03a}.  Up to a constant, $c_n(x;z_c)$
%can be thought of as the probability
%(under a particular critical weighting scheme)
is the probability
that a randomly chosen $n$-step self-avoiding walk
ends at $x$.  Other models to which \cite{HS02} applies include sufficiently
spread-out models of oriented percolation in dimensions $d>4$
\cite{HS03b}, where the corresponding
quantity is the critical two-point function $\tau_n(x;z_c)=\mP((0,0)\ra (x,n))$,
and self-avoiding walks with nearest-neighbour attraction in dimensions $d>4$
\cite{Uelt02}.
More generally, an inductive analysis of lace expansion recursions has been useful
in studying the contact process \cite{HS04z} (extension to continuous time),
self-interacting random walks (such as excited random walk) \cite{HH07} and the
ballistic behavior of 1-dimensional weakly self-avoiding walk \cite{Hofs01}.

As it is stated in \cite{HS02}, the general inductive method is limited to
models with critical dimension $4$.  Thus it does not apply directly to
percolation, which has critical dimension $6$, or to lattice trees, which
have critical dimension $8$.
In this paper, we show that
the method and results of \cite{HS02} are robust
to appropriate changes in various parameters and exponents, so that one can indeed
extend the results to more general
 critical dimensions.

Our extension has been applied already to
prove Gaussian asymptotic behavior for
the two-point function $t_n(x;z_c)$ for sufficiently spread-out lattice trees
in dimensions $d>d_c=8$ in \cite{Holm05,Holm07}.
Up to
a constant, $t_n(x;z_c)$ is the probability (under a particular critical
weighting scheme) that a randomly chosen finite lattice tree contains the
point $x$, with the unique path in the tree from $0$ to $x$ consisting of
exactly $n$ bonds.
The asymptotic behavior of the Fourier transform of the
two-point function provides a first but significant step towards proving convergence
of the finite-dimensional distributions of the associated sequence of measure-valued
processes to those of the canonical measure of super-Brownian motion
\cite{Holm07,HP07}.

A possible future application of our results is to
study the critical two-point function
$\tau_n(x;z_c)$ for sufficiently spread-out percolation in dimensions $d>d_c=6$.
Here, $\tau_n(x;z_c)$ is the probability that $x$ is in the open cluster of the
origin, with the open path of minimum length connecting the origin and $x$
consisting of exactly $n$ bonds, or, alternatively, with the open path of minimum
length connecting the origin and $x$ containing exactly $n$ \emph{pivotal} bonds.

\section{The recursion relation}
The lace expansion typically gives rise to a recursion relation for a \ch{sequence}
$f_n$ depending on parameters $k \in [-\pi,\pi]^d$ and positive $z$.
We may assume that $f_0=1$.  The recursion
relation takes the form
\eqnlab{fkrr}{
        f_{n+1}(k;z) =
        \sum_{m=1}^{n+1} g_m(k;z)f_{n+1-m}(k;z) + e_{n+1}(k;z),
        \quad (n \geq 0),}
with given sequences $g_m(k;z)$ and $e_{n+1}(k;z)$.
The goal is to understand the behaviour of the solution $f_n(k;z)$ of (\ref{fkrr}).

A rough idea of the behaviour we seek to prove can be obtained from the following
(nonrigorous) argument.  Suppose for simplicity that
$D(x)$ is uniformly distributed on a finite box centred at the origin
(so that $\sum_x D(x)=1$),
that $g_1(k;1)=\wh{D}(k)\approx 1-|k|^2\sigma^2/(2d)$, and
that $e_m,g_{m+1}\approx 0$ for $m\ge 1$.
Then we have $f_{n+1}\approx g_1 f_n$, so
$f_n(k) \approx g_1(k)^n\approx \left(1- \frac{|k|^2\sigma^2}{2d}\right)^n$, and thus
\[f_n\left(\frac{k}{\sqrt{\sigma^2 n}};1\right)\approx
\left(1- \frac{|k|^2}{2dn}\right)^n\ra e^{-\frac{|k|^2}{2d}}, \quad \text{ as }n\ra \infty.\]
\ch{
The above argument is, however, overly simplistic, and misses important effects
on the asymptotic behaviour of the solution to \eqref{fkrr} due to
the presence of $e_m(k;z)$ and $g_m(k;z)$.
The inductive method of \cite{HS02} details specific bounds on
$g_m$ and $e_{n+1}$ that ensure that there exists a critical value $z_c$ and
positive constants $A,v$ such that the true asymptotic behaviour is
$f_n\left(\frac{k}{\sqrt{v\sigma^2n}};z_c\right)\ra Ae^{-\frac{|k|^2}{2d}}$.
}
Verification of these bounds has been carried out for
sufficiently spread-out models of
self-avoiding walk \cite{HS03a}, oriented percolation \cite{HS03b} and
the contact process
\cite{HS04z}, by estimating certain Feynman diagrams in dimensions $d>4$.
The required bounds are typically of the form $|h_m(k,z)|\le Cm^{b-\frac{d}{2}}$,
for some functions $h_m$ and exponent $b\ge 0$ that varies from bound to bound.
What turns out to be important in the analysis is that
$\frac{d}{2}=2+\frac{d-4}{2}$ is greater than $2$ when $d>4$.

%In unpublished work \cite{HS02e} the authors note that the analysis of
%\cite{HS02} should be robust enough to permit extension to certain other
%models where the lace expansion is applicable, above $d_c \ne 4$.
%In particular they outline how \cite{HS02} might be adapted to
%analyse lattice trees in dimensions $d>8$.  While deviating
%somewhat in the details, our analysis in this note
%(and its application to lattice trees) is based on the ideas of \cite{HS02e}.

In our analysis we introduce two new parameters $\theta(d)$, $p^*$ and a set $B
\subset [1,p^*]$.
We will discuss the significance of $p^*$ and $B$ following Assumption~D in the
next section.  The most important parameter, $\theta(d)$,
takes the place of $\frac{d}{2}$ in exponents appearing in various bounds.
As in \cite{HS02} we require that $\theta>2$.  In \cite{Holm07},
the result of this note is applied to lattice trees with the choice
$\theta=2+\deight$, with $d>8$.
%We also expect the result to be applicable to other models where the
%lace expansion is used in the analysis above a critical dimension $d_c$.
In general, when the critical dimension is $d_c$, we expect that the correct
parameter value is $\theta=2+\frac{d-d_c}{2}$, e.g.,
we expect that $\theta=2+\frac{d-6}{2}$
is the appropriate choice for percolation.
%Note that in the case $d_c=4$, $2+\frac{d-d_c}{2}=\frac{d}{2}$,
%which is that appearing in \cite{HS02}.
A detailed proof of the results in this note is available in \cite{HHS07x},
% of this consisting of full proofs of the material in \cite{HS02}, adapted to our more general setting,
however, most of the changes to the proof in \cite{HS02} simply
involve replacing $\frac{d}{2}$ in \cite{HS02} with $\theta$ in \cite{HHS07x}.
In this note we  state the new assumptions and results explicitly,
but for the sake of brevity, we present only significant changes in
the proof and refer the reader to \cite{HS02} when the changes are
merely cosmetic.

The remainder of this note is organised as follows.
In Section~\ref{sec:assthm} we state the Assumptions~S, D, E$_{\theta},$
and G$_{\theta}$ on the quantities appearing in the recursion relation,
and the main theorem  to be proved.  In
Section~\ref{sec-ih}, we introduce the induction hypotheses on
$f_n$ that will be used to prove the main theorem.
We then discuss the necessary changes to the advancement of the
induction hypotheses of \cite{HS02}.
% is highly technical and our extension does not require significant alterations to the analysis of .
%therefore briefly discuss the role of $\theta$ in this section and direct the interested reader to \cite{HS02} for the analysis.
Once the induction hypotheses have been advanced,
the main theorem follows without difficulty.

\section{Assumptions and main result}
\label{sec:assthm}

Suppose that for $z>0$ and $k \in [-\pi,\pi]^d$, we have $f_0(k;z) = 1$ and that (\ref{fkrr}) holds for all $n \ge 0$,
%       \eqnlab{fkrec}{
%        f_{n+1}(k;z) =
%        \sum_{m=1}^{n+1} g_m(k;z)f_{n+1-m}(k;z) + e_{n+1}(k;z),
%        \quad \quad (n \geq 0),}
where the functions $g_m$ and $e_m$ are to be regarded as given.
Fix $\theta >2$.

%\subsection{Assumptions S,D,E$_\theta$,G$_\theta$}
The first assumption,
Assumption S, remains unchanged from \cite{HS02}.  It requires that the functions appearing in the recursion relation (\ref{fkrr}) respect the lattice symmetries of reflection and rotation,
and that $f_n$ remains bounded in a weak sense.

\smallskip \noindent
{\bf Assumption S}. For every $n\in\N$ and $z>0$, the mapping
$k\mapsto f_n(k;z)$ is symmetric under replacement of any component
$k_i$ of $k$ by
$-k_i$, and under permutations of the components of $k$.  The same
holds for $e_n(\cdot;z)$ and $g_n(\cdot; z)$.  In addition, for
each $n$, $|f_n(k;z)|$ is bounded uniformly in $k \in
[-\pi,\pi]^d$ and $z$ in a neighbourhood of $1$ (both the bound and the neighbourhood may depend on $n$).

\smallskip
The next assumption, Assumption~D, is only cosmetically changed from \cite{HS02}.
It introduces a probability mass function $D=D_L$ on $\Z^d$ which defines an underlying random walk model and involves a non-negative parameter $L$ which will typically be large.  This serves to spread out the steps of the random walk over a large set.  An example of a family of $D$'s obeying the assumption is taking $D$ uniform on a box of side $2L+1$ centred at the origin.  In particular, Assumption~D implies that $D$ has a finite second moment, and we define
\seqnlab{sigdef}{\sigma^2 \equiv  - \nabla^2 \hat{D}(0)
%&=-\left[\sum_{j=1}^d \frac{\del^2 }{\del k_j^2}\sum_x e^{ik\cdot x}D(x)\right]_{k=0}
=\sum_x |x|^2 D(x),}
%\eqnlab{adef}{a(k) = 1 - \hat{D}(k).}
where $\hat{D}(k) = \sum_{x \in \Z^d} D(x) e^{ik\cdot x}$ is the Fourier transform
of $D$, and $\nabla^2 =\sum_{j=1}^d
\frac{\partial^2}{\partial k_j^2}$ with $k=(k_1,\ldots,k_d)$.

\smallskip \noindent
{\bf Assumption D}.
We assume that
\[f_1(k;z) = z \hat{D}(k) \quad \text{and}\quad e_1(k;z)=0.\]
In particular, this implies that
$g_1(k;z)=z\hat{D}(k)$.  In addition, we also assume:
\\ (i)
$D$ is normalised so that
$\hat{D}(0) =1$, and has $2+2\epsilon$ moments for some $0<\epsilon<\theta-2$, i.e.,
\eqnlab{momentD}{\sum_{x\in \Z^d} |x|^{2+2\epsilon} D(x) <\infty.}
(ii)
There is a constant $C$ such that, for all $L \geq 1$,
\eqnlab{beta, sigmadef}{\|D\|_\infty \leq CL^{-d}\quad \text{and}\quad
\sigma^2  \leq CL^2.}
\\
(iii)
Let $a(k) = 1 - \hat{D}(k)$.  There exist
constants $\eta,c_1,c_2 >0$ such that
\eqnlab{Dbound1}{c_1 L^2 |k|^2 \leq a(k) \leq c_2 L^2 |k|^2 \quad (\|k\|_\infty \leq L^{-1}),}
\eqnlab{Dbound2}{a(k) > \eta  \quad (\|k\|_\infty \geq L^{-1}),}
\eqnlab{Dbound3}{a(k) < 2-\eta \quad (k \in [-\pi,\pi]^d).}

Assumptions~E and~G of \cite{HS02} are adapted to general $\theta>2$ as follows.
The relevant bounds on $f_m$, which {\em a priori}\/ may or may
not be satisfied, are that for some $p^*\ge 1$ and some nonempty $B\subset[1,p^*]$,
we have for every $p \in B$,
\eqnlab{fbdsp}{
\|\hat{D}^2 f_m(\cdot;z)\|_p\leq \frac{K}{L^{\frac{d}{p}} m^{\frac{d}{2p} \wedge \theta}},
    \quad  | f_m(0;z)|\leq K, \quad
    |\nabla^2 f_m(0;z)|\leq K \sigma^2 m,}
for some positive constant $K$.
The bounds in (\ref{fbdsp}) are identical to the ones in \cite[(1.27)]{HS02},
except the first bound, which only appears in \cite{HS02} with $p=1$ and
$\theta=\frac{d}{2}$.
It may be that $B=\{p^*\}$ (i.e. $B$ is a singleton), and then $p=p^*$.
This is the case in \cite{Holm07}, where the choices $p^*=2$ and $B=\{2\}$
are sufficient, as only the $p=2$ case in (\ref{fbdsp}) is required to estimate
the diagrams arising from the lace expansion and verify the assumptions
E$_{\theta}$, G$_{\theta}$ which follow below.
The set $B$ allows for the possibility that
in other applications a larger collection of $\|\cdot\|_p$
norms may be required to verify the assumptions.
Let
\[\beta=\beta(p^*)=L^{-\frac{d}{p^*}}.\]
The parameter $p^*$ serves to make $B$ bounded,
so that $\beta(p^*)$ is small for large $L$.

\smallskip \noindent
{\bf Assumption E$_{\theta}$}. There is an $L_0$, an interval $I \subset
[1-\alpha,1+\alpha]$ with $\alpha \in (0,1)$, and a function $K
\mapsto C_e(K)$, such that if (\ref{fbdsp}) holds for some $K>1$, $L \geq L_0$, $z \in I$ and for all $1 \leq m \leq n$,
then for that $L$ and $z$, and for all $k \in [-\pi,\pi]^d$ and $2
\leq m\leq n+1$, the following bounds hold:
\[|e_m(k;z)|\leq C_e(K) \beta m^{-\theta},
    \quad |e_m(k;z)-e_m(0;z)|\leq
    C_e(K) a(k) \beta m^{-\theta+1}.\]

\smallskip \noindent
{\bf Assumption G$_{\theta}$}. There is an $L_0$, an interval $I \subset
[1-\alpha,1+\alpha]$ with $\alpha \in (0,1)$, and a function $K
\mapsto C_g(K)$, such that if (\ref{fbdsp}) holds for some $K>1$, $L \geq L_0$, $z \in I$ and for all $1 \leq m \leq n$, then for
that $L$ and $z$, and for all $k \in [-\pi,\pi]^d$ and $2 \leq
m\leq n+1$, the following bounds hold:
\[
    |g_m(k;z)|\leq C_g(K)\beta m^{-\theta},
    \quad |\nabla^2 g_m(0;z)|\leq
    C_g(K) \sigma^2 \beta m^{-\theta+1},\]
\[|\partial_z g_m(0;z)|\leq
    C_g(K)\beta m^{-\theta+1},\]
\[|g_m(k;z)-g_m(0;z)- a(k) \sigma^{-2}
    \nabla^2 g_m(0;z)| \leq C_g(K)\beta
    a(k)^{1+\epsilon'}m^{-\theta+1+\epsilon'},\]
with the last bound valid for any $\epsilon' \in [0,\epsilon]$,
with $0<\epsilon<\theta-2$ given by (\ref{momentD}).

\smallskip
Our main result is the following theorem.
(There is a misprint in \cite[Theorem~1.1(a)]{HS02} whose restrictions
should require $\gamma, \delta < \frac{d-4}{2}$ rather than
$\gamma, \delta < \frac{d-4}{4}$;
our assumption $\epsilon< \theta -2$ makes the restriction redundant here.)

\begin{THM}
\label{thm-1p} Let $d>d_c$ and $\theta(d)>2$,
and assume that Assumptions $S$, $D$,
$E_{\theta}$ and $G_{\theta}$ all hold. There exist positive $L_0 = L_0(d,\epsilon)$,
$z_c=z_c(d,L)$, $A=A(d,L)$, and $v = v(d,L)$, such that for $L
\geq L_0$, the following statements hold.
\\
(a)  Fix $\gamma  \in (0,1\wedge \epsilon)$ and $\delta \in (0, (1 \wedge \epsilon) -\gamma)$.
Then
\[f_n\Big(\frac{k}{\sqrt{v \sigma^2 n}};z_c\Big) =A e^{-\frac{|k|^2}{2d}}
        [1 + {\cal O}(|k|^2n^{-\delta})+{\cal O}(n^{-\theta+2})],\]
with the error estimate uniform in $\{k \in \Rd: a(k/\sqrt{v\sigma^2 n}) \leq \gamma n^{-1} \log n \}$.
\\
(b)
\[-\frac{\nabla^2 f_n(0;z_c)}{f_n(0;z_c)}=v \sigma^2 n [1+{\cal O}(\beta n^{-\delta})].\]
(c)  For all $p\ge1$,
\[\|\hat{D}^2 f_n(\cdot;z_c)\|_p\leq \frac{C}{L^{\frac{d}{p}}n^{\frac{d}{2p}\wedge \theta}}.\]
(d) The constants $z_c$, $A$ and $v$ obey
\[1 =  \sum_{m = 1}^\infty g_m(0;z_c), \qquad
    A = \frac{1+\sum_{m=1}^\infty e_m(0;z_c)}{\sum_{m=1}^\infty m g_m(0;z_c)},\qquad
    v =-\frac{\sum_{m=1}^\infty \nabla^2 g_m(0;z_c)}{\sigma^2\sum_{m = 1}^\infty m g_m(0;z_c)}.\]
\end{THM}

As in the proof of \cite[Theorem~1.1]{HS02}, the proof of Theorem~\ref{thm-1p}
establishes the bounds \eqref{fbdsp} for all non-negative integers $m$, with $z$ in
an $m$-dependent interval containing $z_c$.  Consequently, all bounds appearing
in Assumptions E$_\theta$ and G$_\theta$ follow as a corollary, for $z=z_c$ and all
$m$.  Also,
it follows immediately from Theorem~\ref{thm-1p}(d) and the bounds of
Assumptions~E$_\theta$ and G$_\theta$
that
\[z_c=1+ \mc{O}(\beta), \quad A=1+ \mc{O}(\beta), \quad v = 1+ \mc{O}(\beta).\]
Finally, we remark that it is straightforward to extend \cite[Theorem~1.2]{HS02}
for the susceptibility to our present setting, with the
assumption $\theta >2$ replacing $d>4$.  On the other hand,
the proof of the local central
limit theorem \cite[Theorem~1.3]{HS02} does require $\theta = \frac d2$.

\section{Induction hypotheses and their consequences} \label{sec-ih}

\subsection{Induction hypotheses}

Theorem~\ref{thm-1p} is proved
via  induction on $n$, as  in \cite{HS02}.
The induction hypotheses involve a sequence $v_n$, which is
defined exactly as in \cite{HS02} as follows.  We set $v_0=b_0=1$, and for $n \geq
1$ we define
    \[
        b_n = -\frac{1}{\sigma^2}\sum_{m=1}^{n}
    \nabla^2 g_m(0;z)
    ,\quad
    c_n  =  \sum_{m=1}^{n} (m-1) g_m(0;z),\quad
    \label{Delta_n}
    v_n  = \frac{b_n}{1+c_n}.
    \]
The induction hypotheses also involve several constants. Let
$\theta>2$,
and recall from (\ref{momentD}) that $\epsilon<\theta-2$.  We fix
$\gamma, \delta>0$ and $\lambda>2$ according to
\seqnlab{agddef}{0<\gamma<1 \wedge \epsilon, \qquad
0<\delta<(1 \wedge \epsilon)-\gamma, \qquad \theta-\gamma<\lambda&<\theta.}
Here $\lambda$ replaces $\rho+2$ from \cite{HS02}, which is merely a change of notation.

We also introduce constants $K_1, \ldots , K_5$, which are independent of $\beta$.
We define
\eqnlab{K4'def}{
        K_4' = \max \{C_e(cK_4), C_g(cK_4), K_4\},
}
where $c$ is a constant determined in the proof of Lemma~\ref{lem-pibds} below.
To advance the induction, we need to assume that
    \eqnlab{Kcond}{
    K_3 \gg K_1 > K_4' \geq K_4 \gg 1, \quad
    K_2 \geq K_1, 3K_4' , \quad
    K_5 \gg K_4.
    }
Here $a \gg b$ denotes the statement that $a/b$ is sufficiently large.
The amount by which, for instance, $K_3$ must exceed $K_1$ is
independent of $\beta$, but may depend on $p^*$, and is determined during the course of
the advancement of the induction.

Let $z_0=z_1=1$, and define $z_n$ recursively by
    \[
    \label{z_n}
    z_{n+1} = 1-\sum_{m=2}^{n+1}g_m(0;z_n),
    \qquad n \geq 1.
    \]
For $n \geq 1$, we define intervals
    \eqnlab{Indef}{
    I_n = [z_n - K_1\beta n^{-\theta+1}, z_n + K_1\beta n^{-\theta+1}].
    }
In particular this gives $I_1=[1-K_1\beta, 1+K_1\beta]$.

Recall the definition $a(k)=1-\hat{D}(k)$. Our
induction hypotheses are that the following four statements hold
for all $z \in I_n$ and all $1\leq j\leq n$.

\begin{description}
\item[(H1)]
$|z_j - z_{j-1}| \leq K_1 \beta j^{-\theta}$.
\item[(H2)]
$|v_j - v_{j-1}| \leq K_2 \beta j^{-\theta+1}$.
\item[(H3)]
For $k$ such that $a(k) \leq \gamma j^{-1}\log j$, $f_j(k;z)$ can
be written in the form
\[
    f_j(k;z) = \prod_{i=1}^j\left[
    1 -v_i a(k)+ r_i(k)  \right],
\]
with $r_i(k) = r_i(k;z)$ obeying
\[
    |r_i(0)|\leq K_3 \beta i^{-\theta+1},\quad
    |r_i(k)-r_i(0)| \leq K_3 \beta a(k) i^{-\delta}.
\]
\item[(H4)]
For $k$ such that $a(k) > \gamma j^{-1}\log j$, $f_j(k;z)$ obeys
the bounds
\[
    |f_j(k;z)| \leq K_4 a(k)^{-\lambda}j^{-\theta}, \quad
    |f_j(k;z) -f_{j-1}(k;z)| \leq K_5 a(k)^{-\lambda+1} j^{-\theta}.
\]
\end{description}

\medskip \noindent

Note that these four statements are those of \cite{HS02} with the replacement
\eqnlab{eq:replacerho}{\rho+2 \mapsto \lambda}
 in (H4) and the global replacement
\eqnlab{eq:replacetheta}{\frac{d}{2}\mapsto \theta.}
By global replacement we also mean that $\dtwo\mapsto \theta-1$, $\dfour\mapsto \theta-2$, etc. whenever such quantities appear in exponents.

\subsection{Initialisation of the induction}

The verification that the induction hypotheses hold for $n=0$ remains unchanged from the $p=1$ case, up to the replacements (\ref{eq:replacerho}-\ref{eq:replacetheta}).

\subsection{Consequences of induction hypotheses}
\label{sec-prel}
The key result of this section is that the induction hypotheses imply (\ref{fbdsp})
for all $1 \leq m \leq n$, from which the bounds of Assumptions $E_{\theta}$
and $G_{\theta}$ then follow, for $2 \leq m \leq n+1$.

Throughout this note:
    \begin{itemize}
    \item $C$ denotes a strictly positive constant that may depend
    on $d,\gamma,\delta,\lambda$, but {\it not}\/ on the $K_i$, $k$, $n$, and not on
    $\beta$ (provided $\beta$ is sufficiently small, possibly
    depending on the $K_i$). The value of $C$ may change
    from one occurrence to the next.
    \item We frequently assume $\beta \ll 1$ without explicit comment.
    \end{itemize}

Lemmas~\ref{lem-In} and \ref{lem-cA} are proved in \cite{HS02} and the proof in our context requires only the global change (\ref{eq:replacetheta}).

\begin{LEM}
\label{lem-In} Assume (H1) for $1 \leq j \leq n$. Then $I_1
\supset I_2 \supset \cdots \supset I_{n}$.
\end{LEM}

\begin{REM}%[Erratum]
The bound \cite[(2.19)]{HS02} is missing a constant.
Instead of \cite[(2.19)]{HS02} we use
\eqnlab{sbd}{|s_i(k)| \leq K_3 (2+C(K_2+K_3)\beta) \beta a(k) i^{-\delta},}
the only difference being that the constant $2$ appears here instead of a constant
$1$ in \cite[(2.19)]{HS02}.  This does not affect the proof in \cite{HS02}.
To verify (\ref{sbd}), we use the fact that $\R{1-x}\le 1+2x$ for $0\le x \le \hlf$
and note that for small enough $\beta$ it follows from \cite[(2.20)]{HS02} that
\begin{equation*}
\begin{split}
|s_i(k)|&\le \left[1+2K_3\beta\right]\left[(1+|v_i-1|)a(k)r_i(0)+|r_i(k)-r_i(0)|\right]\\
&\le \left[1+2K_3\beta\right]\left[(1+CK_2\beta)a(k)\sm{K_3\beta}{i}{\theta-1}+\sm{K_3\beta a(k)}{i}{\delta}\right]\\
&\le \sm{K_3\beta a(k)}{i}{\delta}[1+2K_3\beta][2+CK_2\beta]\le \sm{K_3\beta a(k)}{i}{\delta}[2+C(K_2+K_3)\beta].
\end{split}
\end{equation*}
Here we have used the bounds of (H2-H3) as well as the fact that $\theta-1>\delta$.
\end{REM}

\begin{LEM}
\label{lem-cA} Let $z\in I_n$ and assume (H2-H3) for $1 \leq
j \leq n$. Then for $k$ with $a(k) \leq \gamma j^{-1}\log j$,
\[
        |f_j(k;z)| \leq e^{CK_3\beta} e^{-(1-C(K_2+K_3)\beta)ja(k)}.\]
\end{LEM}

The middle bound of (\ref{fbdsp}) follows, for $ 1 \leq m \leq n$
and $z \in I_m$, directly
from Lemma~\ref{lem-cA}.  We next state two lemmas which provide the other two bounds of (\ref{fbdsp}).  The first concerns the $\|\cdot\|_p$ norms and contains the most significant changes to \cite{HS02}.  As such we present the full proof of this lemma.

\begin{LEM}
\label{lem-Lpnorm}
Let $z \in I_n$ and assume (H2), (H3) and (H4).
Then for all $1 \leq j \leq n$, and $p \ge 1$,
\[
        \| \hat{D}^2 f_j(\cdot ;z)\|_p \leq \frac{C(1+K_4)}{L^{\frac{d}{p}}j^{\frac{d}{2p}\wedge \theta}},\]
where the constant $C$ may depend on $p,d$.
\end{LEM}

\proof We show that
\[
        \| \hat{D}^2 f_j(\cdot ;z)\|_p^p \leq \frac{C(1+K_4)^p}{L^{d} j^{\frac{d}{2}\wedge \theta p}}.\]
For $j=1$ the result holds since $|f_1(k)|=|z\wh{D}(k)|\le z\le 2$, and,
since $p\geq 1$, it therefore follows from (\ref{beta,  sigmadef}) and the Parseval
relation that
$ \| \hat{D}^2 f_1(\cdot ;z)\|_p^p \leq 2^p \|\hat D ^{2p}\|_1
\leq 2^p\|\hat D ^2\|_1 = 2^p\|D\|_2^2 \leq 2^pCL^{-d}$.  We may therefore assume that $j\ge 2$ where needed in what follows, so that in particular $\log j \ge \log 2$.

Fix $z \in I_n$ and $1 \leq j \leq n$, and define
\begin{equation*}
\begin{split}
    R_{1} & = \{k \in
    [-\pi,\pi]^{d}: a(k) \leq \gamma j^{-1}\log j , \;
    \|k\|_\infty \leq L^{-1} \} ,\nonumber \\
    R_{2} & = \{k \in
    [-\pi,\pi]^{d}: a(k) \leq \gamma j^{-1}\log j , \;
    \|k\|_\infty > L^{-1} \} ,\nonumber \\
    \label{Rjdef}
    R_{3} & = \{k \in
    [-\pi,\pi]^{d}: a(k) > \gamma j^{-1}\log j, \;
    \|k\|_\infty \leq L^{-1} \},\nn\\
    R_4 & = \{k \in
    [-\pi,\pi]^{d}: a(k) > \gamma j^{-1}\log j, \;
    \|k\|_\infty > L^{-1} \}.
    \end{split}
\end{equation*}
The set $R_2$ is empty if $j$ is sufficiently large.  Then
    \[
    \|\hat{D}^2f_{j}\|_p^p = \sum_{i=1}^4\int_{R_i}\left(\hat{D}(k)^2|f_{j}(k)|\right)^p
    \frac{d^dk}{(2\pi)^d}.
    \]
We will treat each of the four terms on the right side separately.

On $R_1$, we use (\ref{Dbound1}) in conjunction with
Lemma~\ref{lem-cA} and the fact that \ch{$\hat{D}(k)^2\leq 1$}, to obtain for all $p>0$,
\begin{equation*}
\begin{split}
    \int_{R_1}\left(\hat{D}(k)^2|f_{j}(k)|\right)^p \frac{d^dk}{(2\pi)^d}
    &\leq \int_{R_1} Ce^{-cpj(L|k|)^2} \frac{d^dk}{(2\pi)^d}\\
&\le  \int_{{\mathbb R}^d} Ce^{-cpj(L|k|)^2}dk
    \leq \frac{C}{L^d (pj)^{d/2}}\le \frac{C}{L^d j^{d/2}}.
    \end{split}
\end{equation*}
Here we have used the substitution $k'_i=Lk_i\sqrt{pj}$.
On $R_2$, we use Lemma~\ref{lem-cA} and (\ref{Dbound2}) to
conclude that for all $p>0$, there is an $\alpha(p) > 1$ such that
    \eqn{
    \int_{R_2} \left(\hat{D}(k)^2|f_{j}(k)|\right)^p \frac{d^dk}{(2\pi)^d}
    \leq C\int_{R_2} \alpha^{-j} \frac{d^dk}{(2\pi)^d}
    = C\alpha^{-j} |R_2|,}
where $|R_2|$ denotes the volume of $R_2$.  This volume is maximal
when $j=3$, so that

    \seqn{
    |R_2| &\leq \Big|\{ k : a(k)\le \textstyle\frac{\gamma \log 3}{3} \}\Big|
    \leq \Big|\{k: \hat{D}(k) \geq 1-\textstyle\frac{\gamma \log 3}{3} \}\Big| \leq \Big(\textstyle\frac{1}{1-\frac{\gamma \log 3}{3}}\Big)^2 \|\hat{D}^2\|_1
    \leq \Big(\textstyle\frac{1}{1-\frac{\gamma \log 3}{3}}\Big)^2 CL^{-d},}
using (\ref{beta, sigmadef}) in the last step.
Therefore $\alpha^{-j}|R_2| \leq C L^{-d} j^{-d/2}$ since $\alpha^{-j}j^{\frac{d}{2}}\le C(\alpha,d)$ for every $j$, and
    \eqn{
    \int_{R_2} \left(\hat{D}(k)^2|f_{j}(k)|\right)^p \frac{d^dk}{(2\pi)^d}
    \leq CL^{-d} j^{-d/2}.
    }

On $R_3$ and $R_4$, we use (H4).  As a result, the contribution
from these two regions is bounded above by
    \eqn{
    \left(\frac{K_4 }{j^\theta}\right)^p \sum_{i=3}^4\int_{R_i}
        \frac{\hat{D}(k)^{2p}}{a(k)^{\lambda p}} \frac{d^dk}{(2\pi)^d}.}
We first consider $R_3$, where we apply $\hat{D}(k)^2\leq 1$.
Recall that we can restrict our attention to $j\ge 2$.
From (\ref{Dbound1}), $k \in R_3$ implies that $L^2|k|^2>Cj^{-1}\log j$, and
we have the upper bound
\begin{equation}
    \label{largekext}
    \frac{CK_4^p}{j^{\theta p} L^{2\lambda p}}
    \int_{R_3} \frac{1}{|k|^{2\lambda p}} d^dk
\le \frac{CK_4^p}{j^{\theta p} L^{2\lambda p}}
\int_{\sqrt{\frac{C\log j}{L^2j}}}^{\frac{C}{L}}r^{d-1-2\lambda p}dr.
\end{equation}
%Since $\log 1 =0$, this integral will not be finite if both $j=1$
%and $p \ge \frac{d}{2\lambda }$, but
For $d>2\lambda p$, we have an upper bound on (\ref{largekext}) of
\eqnlab{eq:dgreater}{\frac{CK_4^p}{j^{\theta p} L^{2\lambda p}}
\int_0^{\frac{C}{L}}r^{d-1-2\lambda p}dr
\le \frac{CK_4^p}{j^{\theta p} L^{2\lambda p}}
\left(\frac{C}{L}\right)^{d-2\lambda p}\le \frac{CK_4^p}{j^{\theta p}L^d}.}
For $d=2\lambda p$, (\ref{largekext}) is
\seqnlab{eq:dequal}{\frac{CK_4^p}{j^{\theta p} L^{2\lambda p}}
\int_{\sqrt{\frac{C\log j}{L^2j}}}^{\frac{C}{L}}\R{r}dr
\le \frac{CK_4^p}{j^{\theta p} L^{2\lambda p}}
\log \left(\frac{C\sqrt{L^2j}}{L\sqrt{\log j}}\right)
=\frac{CK_4^p}{j^{\theta p} L^{2\lambda p}}\log\left(\frac{Cj}{\log j}\right),}
and $\theta p=\frac{\theta d}{2\lambda}>\frac{d}{2}$
%Write $\frac{CK_4^p}{j^{\theta p}L^d}=\frac{CK_4^p}{j^{\frac{d}{2}+(\frac{\theta}{\lambda}-1)}L^d}$ and
since $\lambda <\theta$.   This gives an upper bound
in this case of $CK_4^pj^{-\frac{d}{2}}L^{-d}$.
Lastly, for $d<2\lambda p$, since $\lambda<\theta$, (\ref{largekext}) is bounded,
as required, by
\eqnlab{eq:dless}{\frac{CK_4^p}{j^{\theta p} L^{2\lambda p}}
\int_{\sqrt{\frac{C\log j}{CL^2j}}}^{\infty}r^{d-1-2\lambda p}dr
\le \frac{CK_4^p}{j^{\theta p} L^{2\lambda p}}
\left(\frac{CL^2j}{\log j}\right)^{\frac{2\lambda p-d}{2}}
\le \frac{CK_4^p}{j^{\frac{d}{2}}L^d}.}

On $R_4$, we use (\ref{beta, sigmadef}),
$p \geq 1$, \ch{$\hat{D}(k)^2\leq 1$}, and (\ref{Dbound2}) to obtain the bound
\eqn{
\frac{CK_4^p}{j^{\theta p}}
\int_{[-\pi,\pi]^d}  \hat{D}(k)^{2p}\frac{d^dk}{(2\pi)^d}
\leq \frac{CK_4^p}{j^{\theta p}}
\int_{[-\pi,\pi]^d}  \hat{D}(k)^{2}\frac{d^dk}{(2\pi)^d}
\le \frac{CK_4^p}{j^{\theta p}L^d}.}
%Since $K_4^p \le (1+K_4)^p$,
This completes the proof.
\qed

\begin{LEM}
\label{lem-fder}
Let $z \in I_n$ and assume (H2) and (H3).  Then, for
$1 \leq j \leq n$,
\eqn{
        | \nabla^2 f_j(0 ;z) | \leq (1+C(K_2 + K_3) \beta ) \sigma^2 j.
}
\end{LEM}
The proof is identical to \cite{HS02}.  We merely point out one
inconsequential correction to the first line of \cite[(2.35)]{HS02}:
a constant $2$ is missing and  it should read
\eqnlab{1.2b3}{
     \nabla^2 s_i(0)
= 2  \sum_{l=1}^d \lim_{t \rightarrow 0}\frac{s_i(te_l)-s_i(0)}{t^2}  .}

The next lemma, whose proof proceeds exactly as in \cite{HS02} with $\frac{d}{2}$ replaced by $\theta$, is the key to advancing the induction, as it
provides bounds for $e_{n+1}$ and $g_{n+1}$.
Recall that $K_4'$ was defined in \eqref{K4'def}.

\begin{LEM}
\label{lem-pibds}
 Let $z\in I_{n}$, and assume (H2),
(H3) and (H4). For $k \in [-\pi,\pi]^d$, $2 \leq j \leq n+1$, and
$\epsilon' \in [0,\epsilon]$, the following hold:
\vspace{-2mm}
\begin{tabbing}
(iii) \= \kill (i) \>  $|g_j(k;z)|\leq  K_4' \beta j^{-\theta}$,
\\
(ii) \> $|\nabla^2 g_j(0;z)|\leq   K_4'  \sigma^2 \beta
j^{-\theta+1}$,
\\
(iii) \> $|\partial_z g_j(0;z)|\leq  K_4' \beta j^{-\theta+1},$
\\
(iv) \>  $|g_j(k;z)-g_j(0;z)- a(k) \sigma^{-2}\nabla^2 g_j(0;z)|
\leq  K_4' \beta a(k)^{1+\epsilon'}j^{-\theta+1+\epsilon'},$
\\
(v)   \> $|e_j(k;z)|\leq  K_4' \beta j^{-\theta}$,
\\
(vi)  \> $|e_j(k;z)-e_j(0;z)|\leq  K_4' a(k) \beta
j^{-\theta+1}.$
\end{tabbing}
\end{LEM}

\section{The induction advanced}
\label{sec-adv}

The advancement of the induction is carried out as in \cite{HS02} with a
few minor changes corresponding to the global replacement (\ref{eq:replacetheta}),
and also (\ref{eq:replacerho}) for (H4).
Full details can be found in \cite{HHS07x}, and here we only point out the main
places where changes are required.

In adapting \cite[(3.2)]{HS02}, we use the fact that
$\sum_{m=2}^{\infty}m^{-\theta+1}<\infty$,
since $\theta>2$, and in adapting \cite[(3.26)]{HS02}, we use
$\sum_{j=n+2-m}^n j^{-\theta+1}\le C(n+2-m)^{-\theta+2}$.
For \cite[(3.40)]{HS02}, we apply
$\epsilon'\le \epsilon < \theta-2$ to conclude that
$\sum_{m=2}^{\infty}m^{-\theta +1 +\epsilon'}<\infty$.
To adapt \cite[(3.43)]{HS02},
we use the fact that $\delta + \gamma < 1 \wedge (\theta-2)$, by (\ref{agddef}),
to conclude that there exists a $q>1$
sufficiently close to $1$ so that
\eqn{(n+1)^{-\delta} \geq
    (n+1)^{\gamma q - 1} \log (n+1)
        \times \begin{cases} (n+1)^{0 \vee (3-\theta)}, & (\theta \neq 3)
        \\ \log(n+1), & (\theta=3).  \end{cases}}
Other similar bounds required to verify (H3)
(corresponding to \cite[(3.50)--(3.51)]{HS02} and \cite[(3.58)]{HS02} for example)
also follow from $\delta + \gamma < 1 \wedge (\theta-2)$.
For (H4), using the fact that $\gamma+\lambda-\theta>0$, there exists
$q'$ close to $1$ so that for $a(k)\leq \gamma n^{-1} \log n$,
\eqn{\frac{C}{n^{\theta}}\frac{n^{\lambda}}{n^{q'\gamma +\lambda -\theta}}
        \leq \frac{C}{n^{\theta}a(k)^{\lambda}}.}
This corresponds to \cite[(3.62)]{HS02},
and is used to advance the first and second bounds of (H4).
%In addition we use the fact that $\lambda>2$ so that $a(k)^{\lambda-2}\le C$
%(recall that $a(k)\le 2$ from (\ref{Dbound3})) to get $a(k)^{-1}\le Ca(k)^{1-\lambda}$.

\smallskip
Once the induction has been advanced, the proof of Theorem \ref{thm-1p} is
then completed
exactly as in \cite{HS02}, with the global replacement (\ref{eq:replacetheta}).
%%%%%%%%%%%%%%%%%%%%%%%%%%%

\section*{Acknowledgements}
A version of this work appeared in the PhD thesis \cite{Holm05}.  The work of RvdH and MH was supported in part by Netherlands Organisation for Scientific Research (NWO).  The work of GS was supported in part by NSERC of Canada.

%\bibliography{../bibdef/bib}
%\bibliographystyle{plain}

\end{document}